%% file: dec2019.tex
\documentclass[12pt]{article}
\usepackage[english]{babel}
\usepackage{lineno} 
\usepackage{graphicx,multicol}
\usepackage{epic,eepic,epsfig}
\usepackage{caption}

\usepackage{amsmath,amsfonts,amsthm,amssymb}
\usepackage{pstricks,pstricks-add,pst-math,pst-xkey}

\usepackage{algorithm, algorithmic}

\usepackage{graphicx,pstricks,pst-node,amssymb}

\usepackage{tikz}

\pgfdeclarelayer{edgelayer}
\pgfdeclarelayer{nodelayer}
\pgfsetlayers{edgelayer,nodelayer,main}

\tikzstyle{none}=[inner sep=0pt]
\definecolor{hexcolor0xf81e1c}{rgb}{0.973,0.118,0.110}
\definecolor{hexcolor0x3c00ff}{rgb}{0.235,0.000,1.000}

\tikzstyle{whitevertex}=[circle,fill=white,draw=black, scale = 0.5]
\tikzstyle{redvertex}=[circle,fill=hexcolor0xf81e1c,draw=black, scale = 0.5]
\tikzstyle{bluevertex}=[circle,fill=hexcolor0x3c00ff,draw=black, scale = 0.5]
\tikzstyle{greenvertex}=[circle,fill=green,draw=black, scale=0.5]
\tikzstyle{purplevertex}=[circle,fill=magenta,draw=black, scale=0.5]
\tikzstyle{grayvertex}=[rectangle,fill=black,draw=gray, scale=0.5]
\tikzstyle{blackvertex}=[circle,fill=black,draw=black, scale=0.5]

\tikzstyle{textbox}=[rectangle,fill=none,draw=none]
\tikzstyle{box}=[rectangle,fill=none,draw=black]

\tikzstyle{arc}=[black, ->]
\tikzstyle{grayarc}=[gray, ->]
\tikzstyle{bluearc}=[blue, ->]
\tikzstyle{grayedge}=[draw=gray]
\tikzstyle{blueedge}=[draw=blue]
\tikzstyle{rededge}=[draw=red]
\tikzstyle{edge}=[draw=black]

\tikzstyle{vertex}=[circle, ,fill=white,draw=black, scale=0.5]

\tikzstyle{10circle}=[circle, scale=10.0,draw=black]
\tikzstyle{10oval}=[ellipse, scale=10.0,draw=black]

\usepackage{amsmath,amsfonts,amsthm,amssymb}
\usepackage{pstricks,pstricks-add,pst-math,pst-xkey}

\usepackage{algorithm, algorithmic}

\usepackage{graphicx,pstricks,pst-node,amssymb}

\usepackage{algorithm, algorithmic}

\begin{document}

\newtheorem{thm}{\hspace{5mm}Theorem}[section]
\newtheorem{prop}[thm]{\hspace{5mm}Proposition}
\newtheorem{defn}[thm]{\hspace{5mm}Definition}
\newtheorem{lem}[thm]{\hspace{5mm}Lemma}
\newtheorem{cor}[thm]{\hspace{5mm}Corollary}
\newtheorem{conj}[thm]{\hspace{5mm}Conjecture}
\newtheorem{prob}[thm]{\hspace{5mm}Problem}
\newtheorem{quest}[thm]{\hspace{5mm}Question}
\newtheorem{alg}[thm]{\hspace{5mm}Algorithm}
\newtheorem{sub}[thm]{\hspace{5mm}Algorithm}
\newcommand{\induce}[2]{\mbox{$ #1 \langle #2 \rangle$}}
\newcommand{\2}{\vspace{2mm}}
\newcommand{\dom}{\mbox{$\rightarrow$}}
\newcommand{\ndom}{\mbox{$\not\rightarrow$}}
\newcommand{\compdom}{\mbox{$\Rightarrow$}}
\newcommand{\cdom}{\compdom}
\newcommand{\sdom}{\mbox{$\Rightarrow$}}
\newcommand{\lsd}{locally semicomplete digraph}
\newcommand{\lt}{local tournament}
\newcommand{\la}{\langle}
\newcommand{\ra}{\rangle}
\newcommand{\pf}{{\hspace{1mm}\bf Proof: }}
\newtheorem{claim}{Claim}
\newcommand{\beq}{\begin{equation}}
\newcommand{\eeq}{\end{equation}}
\newcommand{\<}[1]{\left\langle{#1}\right\rangle}

\newcommand{\set}[1]{\left\{#1 \right\}}
\newcommand{\propset}[2]{\left\{ #1 \, \left| \, #2 \right.\right\}}
\newcommand{\Z}{\mathbb{Z}}
\newcommand{\R}{\mathbb{R}}
\renewcommand{\cal}[1]{\mathcal{#1}}
\newcommand{\n}{1,2,\dots,n}
\newcommand{\no}{0,1,\dots,n-1}
\newcommand{\seq}[2]{#1_1,#1_2,\dots,#1_{#2}}
\newcommand{\seqo}[2]{#1_0,#1_1,\dots,#1_{#2-1}}
\newcommand{\abs}[1]{\left| #1 \right|}
\newcommand{\ceil}[1]{\left\lceil #1 \right\rceil}
\newcommand{\modular}[1]{\,(\text{mod } #1)}
\newcommand{\vect}[1]{\left(\substack{#1}\right)}
\newcommand{\rel}[2]{\sim_{_{#1,#2}}}
\newcommand{\tdots}{\cdot \ldots \cdot}

\newcommand{\func}[5]{\begin{alignat*}{2}#1:\,#2&\rightarrow #3\\#4&\mapsto#5.\end{alignat*}}

\newcommand{\conditional}[2]{\ifthenelse{\boolean{#1}}{#2}{}}
\newboolean{proofs}
\newboolean{remarks}
\newboolean{personal}
\newboolean{figures}
\setboolean{proofs}{true}
\setboolean{remarks}{true}
\setboolean{personal}{true}
\setboolean{figures}{true}

\bibliographystyle{amsplain}
\addcontentsline{toc}{chapter}{Bibliography}


\title{$(k,l)$-Colourings and Ferrers Diagram Representations of Cographs}

\author{Dennis A. Epple and Jing Huang\thanks{Department of Mathematics and Statistics,
      University of Victoria, Victoria, B.C., Canada V8W 2Y2; huangj@uvic.ca}
      }
\date{}

\maketitle

\begin{abstract}
For a pair of natural numbers $k, l$, a $(k,l)$-colouring of a graph $G$ is 
a partition of the vertex set of $G$ into (possibly empty) sets $\seq{S}{k}$, 
$\seq{C}{l}$ such that each set $S_i$ is an independent set and each set $C_j$ 
induces a clique in $G$. The $(k,l)$-colouring problem, which is NP-complete
in general, has been studied for special graph classes such as chordal graphs,
cographs and line graphs. 
Let $\hat{\kappa}(G) = (\kappa_0(G),\kappa_1(G),\dots,\kappa_{\theta(G)-1}(G))$
and $\hat{\lambda}(G) = (\lambda_0(G),\lambda_1(G),\dots,\lambda_{\chi(G)-1}(G))$
where $\kappa_l(G)$ (respectively, $\lambda_k(G)$) is the minimum $k$
(respectively, $l$) such that $G$ has a $(k,l)$-colouring. 
We prove that $\hat{\kappa}(G)$ and $\hat{\lambda}(G)$ are a pair of 
conjugate sequences for every graph $G$ and when $G$ is a cograph, 
the number of vertices in $G$ is equal to the sum of the entries in 
$\hat{\kappa}(G)$ or in $\hat{\lambda}(G)$. Using the decomposition property of
cographs we show that every cograph can be represented by Ferrers diagram.
We devise algorithms which compute $\hat{\kappa}(G)$ for cographs $G$ and find
an induced subgraph in $G$ that can be used to certify the 
non-$(k,l)$-colourability of $G$.  
\end{abstract}

\section{Introduction}

Let $G$ be a graph and $k, l \geq 0$ be natural numbers. 
A {\em $(k,l)$-colouring} of $G$ is a partition of the vertex set of $G$ into 
(possibly empty) sets $\seq{S}{k}$, $\seq{C}{l}$ such that each $S_i$ is 
an independent set and each $C_j$ induces a clique in $G$. 
The concept of $(k,l)$-colourings encompasses the classical colouring and clique 
covering of graphs; indeed, a $(k,0)$-colouring is just a $k$-colouring 
and a $(0,l)$-colouring is a partition of $G$ into at most $l$ cliques. 
A graph is {\it $(k,l)$-colourable} if it has a $(k,l)$-colouring. 
Thus bipartite graphs are exactly the $(2,0)$-colourable graphs and 
split graphs are precisely the $(1,1)$-colourable graphs \cite{77FoHa}.

The {\em bichromatic number} $\chi^b(G)$ of $G$ is the least integer $r$ such
that, for all $k, l$ with $k+l = r$, $G$ is $(k,l)$-colourable.
The notion of bichromatic number arose in the study of extremal graphs, motivated
by classical results of Tur\'an \cite{turan} and of Erd\H{o}s, Stone and Simonovits
(cf. \cite{46ErSt}). This parameter has been studied by Pr\"omel and Steger
\cite{93PrSt} under the name of {\em $\tau$-parameter}, by Bollob\'as and Thomason
\cite{97BoTh} under the name of {\em colouring number}, and by Axenovich, {K\'ezdy} and Martin \cite{08AxKeMa} under the name of {\em binary chromatic number}.
The parameter is tied to the speed of hereditary properties and edit distance,
cf. \cite{08AlSt,08AxKeMa,97BoTh}. A counterpart of the bichromatic number 
$\chi^b(G)$ is the notion of the {\em cochromatic number} $\chi^c(G)$, which is 
the least integer $r$ such that $G$ is $(k,l)$-colourable for some $k, l$ with 
$k+l = r$, cf. \cite{77LeSt}. 

It is not surprising that computing the bichromatic number of a graph is an NP-hard
problem, cf.  \cite{10EpHu}. Brandst\"adt \cite{96Br,96Br_a} 
proved that the problem of deciding whether a graph is $(k,l)$-colourable is 
NP-complete for fixed $k, l$ with $k \geq 3$ or $l \geq 3$ and polynomial time 
solvable otherwise. A graph is {\em chordal} if it does not contain 
an induced $C_k$ for each $k \geq 4$ and is a {\em cograph} if it does not 
contain an induced $P_4$. It is proved in \cite{04HeKlNoPr} that a chordal graph 
is $(k,l)$-colourable if and only if it does not contain $(l+1)K_{k+1}$, 
the disjoint union of $l+1$ copies of $K_{k+1}$ (see definition below). 
The $(k,l)$-colourability of cographs and line-graphs have been studied in 
\cite{05DeEkWe1,05DeEkWe2,epple,10EpHu,07FeHeHo}.

Every graph $G$ satisfies $\chi^b(G) \leq \chi(G) + \theta(G) - 1$ where $\chi(G)$ 
and $\theta(G)$ are the chromatic number and the clique covering number of $G$ 
respectively, cf. \cite{08AlSt,93PrSt}. 
Graphs which satisfy this inequality with equality have been characterized in 
\cite{10EpHu}. It turns out that all these graphs are cographs. To describe   
the characterization, we recall the recursive definition of cographs. 
The {\em disjoint union} of graphs $G$ and $H$, denoted by
$G + H$, is the graph with vertex set $V(G) \cup V(H)$ and edge set 
$E(G) \cup E(H)$. 

Let $\cal C$ be the set of graphs defined recursively as follows:
\begin{itemize}
\item $K_1 \in \cal C$;
\item if $G \in \cal C$, then $\overline{G} \in \cal C$;
\item if $G, H \in \cal C$, then $G+H \in \cal C$.
\end{itemize}

\begin{thm} \cite{81CoLeSt} \label{corneil}
The following statements are equivalent for a graph $G$:
\begin{enumerate}
\item $G \in \cal C$;
\item $G$ is a cograph (i.e., $G$ does not contain an induced $P_4$);
\item for every induced subgraph $H \neq K_1$ of $G$, either $H$ or $\overline{H}$
      is disconnected.
\qed
\end{enumerate}
\end{thm}   

Let $\cal B$ be a set of graphs constructed recursively as follows:
\begin{itemize}
\item $K_1 \in \cal B$;
\item if $G \in \cal B$, then $\overline{G} \in \cal B$;
\item if $G, H \in \cal B$ with $\chi(G) = \chi(H)$, then $G+H \in \cal B$.
\end{itemize}

The {\em join} of $G$ and $H$, denoted by $G \vee H$, is the graph with 
vertex set $V(G) \cup V(H)$ and edge set 
$E(G) \cup E(H) \cup \{uv:\ u \in V(G),\ v \in V(H)\}$.
It is easy to see that if $G, H \in \cal C$ then $G \vee H \in \cal C$
and if $G, H \in \cal B$ with $\theta(G) = \theta(H)$ then $G \vee H \in \cal B$.

Clearly, ${\cal B} \subseteq {\cal C}$. We call the graphs in $\cal B$ 
{\em box cographs}. It is proved in \cite{10EpHu} that the box cographs are 
exactly the graphs $G$ which satisfy the inequality 
$\chi^b(G) \leq \chi(G) + \theta(G) - 1$ with equality. 

\begin{thm} \cite{10EpHu} \label{notejgt}
A graph $G$ satisfies $\chi^b(G) = \chi(G) + \theta(G) - 1$ if and only if it is
a box cograph. 
\qed
\end{thm}

A box cograph is of {\em dimension $k$ times $l$} if it has chromatic number $k$ 
and clique covering number $l$. The graph $lK_k$ is a box cograph of dimension 
$k$ times $l$. No box cograph of dimension $k+1$ times $l+1$ is 
$(k,l)$-colourable. The following theorem asserts that they are exactly the 
forbidden (induced) subgraphs for $(k,l)$-colourable cographs. 
An equivalent statement of the theorem is proved in \cite{07FeHeHo}.

\begin{thm} \cite{10EpHu} \label{notejgt+}
A cograph is $(k,l)$-colourable if and only if does not contain a box cograph of
dimension $k+1$ times $l+1$ as an induced subgraph.
\qed
\end{thm}

\begin{cor} 
For any cograph $G$,
\[\chi^b(G) = \mbox{max} \{k+l:\ G\ \mbox{contains\ a\ box\ cograph\ of\
        dimension}\ k\ \mbox{times}\ l\}\ - 1\]
\qed
\end{cor}

Let $G$ be a graph and $k, l$ be natural numbers. We use $\kappa_l(G)$ to
denote the minimum $k$ for which $G$ is $(k,l)$-colourable and use
$\lambda_k(G)$ to denote the minimum $l$ for which $G$ is $(k,l)$-colourable.
Write
 \[\hat{\kappa}(G) = (\kappa_0(G),\kappa_1(G),\dots,\kappa_{\theta(G)-1}(G))\]
and
 \[\hat{\lambda}(G) = (\lambda_0(G),\lambda_1(G),\dots,\lambda_{\chi(G)-1}(G)).\]
The knowledge of the values of $\hat{\kappa}(G)$ or in $\hat{\lambda}(G)$ can be 
directly used to determine whether $G$ is $(k,l)$-colourable and hence to compute
the bichromatic number $\chi^b(G)$. Indeed, a graph $G$ is $(k,l)$-colourable 
if and only if $\kappa_l(G) \leq k$ or equivalently, if and only if 
$\lambda_k(G) \leq l$. 

In this paper, we show that $\hat{\kappa}(G)$ and $\hat{\lambda}(G)$ are a pair of 
conjugate sequences for every graph $G$, that is, the Ferrers diagrams of 
$\hat{\kappa}(G)$ and $\hat{\lambda}(G)$ are conjugate to each other. 
We prove that, when $G$ is a cograph, the number of vertices in $G$ is equal to 
the sum of the entries in either of the sequences. This is not true in general.
Using the decomposition property of cographs, we show that every cograph can be
drawn in a shape similar to a Ferrers diagram for a sequence of numbers.
We devise efficient bottom-up and top-down algorithms on cotrees of cographs. 
The bottom-up algorithm calculates the sequence $\hat{\kappa}(G)$ and 
the top-down algorithm finds a box cograph of specified dimension which certifies 
the non-$(k,l)$-colourability of the input graph $G$.  

Algorithms for the $(k,l)$-colourability of cographs have been studied by Gimbel,
Kratsch, and Stewart \cite{94GiKrSt} and by Demange, Ekim and de Werra  
\cite{05DeEkWe1}. In \cite{94GiKrSt}, an algorithm for the computation of 
the cochromatic number of a cograph using its cotree was presented. 
This algorithm, which runs in time $O(n^2)$, implicitly uses $(k,l)$-colourings
of cographs. Demange, Ekim and de Werra \cite{05DeEkWe1} gave two different 
algorithms concerning the $(k,l)$-colouring of cographs. The first, which also 
uses cotrees, calculates a maximum $(k,l)$-colourable induced subgraph of a cograph
(thereby also checking the $(k,l)$-colourability of the cograph itself) in time 
$O((k^3l+kl^3)n)$. For the purpose of the second algorithm, it is shown that if 
$G$ is a $(k,l)$-colourable cograph (with $l \geq 1$) and $C$ a maximum clique of 
$G$, then $G-C$ is $(k,l-1)$-colourable. Using this, the algorithm finds 
a $(k,l)$-colouring of a cograph (if one exists) by successively removing 
$l$ maximum cliques and finding a $k$-colouring of the remaining graph. 
The algorithm runs in time $O(n(m+n))$ where $m$ is the number of edges of 
the graph. An adaptation of this idea is presented which calculates 
the cochromatic number in time $O(n^{3/2})$. Our algorithms for computing 
$\hat{\kappa}(G)$ and for finding certificates for the non-$(k,l)$-colourability
of cographs $G$ matches the same complexity $O(n^2)$ as the algorithms in 
\cite{05DeEkWe1,94GiKrSt}. In fact, the algorithm for 
$\hat{\kappa}(G)$ can be implemented to run in time $O(n \log n)$.

\section{Ferrers diagrams}
\label{Young}

Since $\kappa_0(G) \geq \kappa_1(G) \geq \dots \geq \kappa_{\theta(G)-1}(G)$ and 
$\lambda_0(G) \geq \lambda_1(G) \geq \dots \geq \lambda_{\chi(G)-1}(G)$, 
$\hat{\kappa}$ and $\hat{\lambda}$ are both monotonically non-increasing sequences, we can represent each of them by a {\it Ferrers diagram}.
For example, the graph in Figure~\ref{g1} has $\hat{\kappa} = (3,3,1)$ and 
$\hat{\lambda} = (3,2,2)$, which are represented by Ferrers diagrams in 
Figure~\ref{sqs_g1}. Note that these two sequences are conjugate to each other,
that is, by reflecting any of the diagrams along the main diagonal 
we obtain the other diagram. We show below this is the case for every graph. 

\vspace{-4mm}
\begin{center}
\begin{figure}[htb]
\center
\begin{tikzpicture}[>=latex]
\begin{pgfonlayer}{nodelayer}
                \node [style=blackvertex] (1) at (-1,.5){};
                \node [style=blackvertex] (2) at (1,.5){};
                \node [style=blackvertex] (3) at (0,-.1){};
                \node [style=blackvertex] (4) at (0,-1){};
                \node [style=blackvertex] (5) at (0,-2){};
                \node [style=blackvertex] (6) at (-1,-1.5){};
                \node [style=blackvertex] (7) at (1,-1.5){};
                \draw[- ] (1) -- (2) -- (3) -- (6) -- (5) -- (7) -- (3) -- (1);
                \draw[- ] (4) -- (5);
                \draw[- ] (6) -- (4) -- (7);
\end{pgfonlayer}
\end{tikzpicture}
\caption{\label{g1}A graph}
\end{figure}
\end{center}

\vspace{-4mm}
\begin{center}
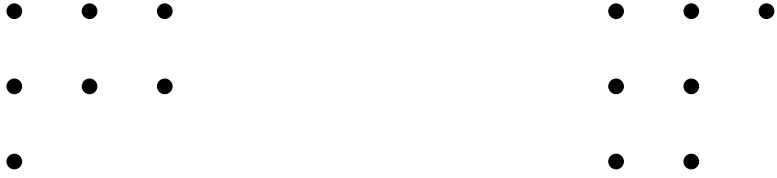
\begin{figure}[htb]
\center
\begin{tikzpicture}[>=latex]
\begin{pgfonlayer}{nodelayer}
                \node [style=blackvertex] (1) at (0,0){};
                \node [style=blackvertex] (2) at (1,0){};
                \node [style=blackvertex] (3) at (2,0){};
                \node [style=blackvertex] (4) at (0,-1){};
                \node [style=blackvertex] (5) at (1,-1){};
                \node [style=blackvertex] (6) at (2,-1){};
                \node [style=blackvertex] (7) at (0,-2){};

                \node [style=blackvertex] (8) at (8,0){};
                \node [style=blackvertex] (9) at (9,0){};
                \node [style=blackvertex] (10) at (10,0){};
                \node [style=blackvertex] (11) at (8,-1){};
                \node [style=blackvertex] (12) at (9,-1){};
                \node [style=blackvertex] (13) at (8,-2){};
                \node [style=blackvertex] (14) at (9,-2){};
\end{pgfonlayer}
\end{tikzpicture}
\caption{\label{sqs_g1}Ferrers diagrams of $\hat{\kappa}$ and 
           $\hat{\lambda}$ for the graph in Figure \ref{g1}}
\end{figure}
\end{center}

\begin{thm} \label{conjugate}
For any graph $G$, $\hat{\kappa}(G)$ and $\hat{\lambda}(G)$ 
are a pair of conjugate sequences.
\end{thm}
\pf For each $j$ with $0 \leq j \leq \chi(G) - 1$, let 
\[i^*_j = \mbox{max}\{i: \ \kappa_i(G) \geq j+1\}. \]
To prove that $\hat{\lambda}(G)$ is the conjugate sequence of $\hat{\kappa}(G)$,
it suffices to show that $\lambda_j(G) = i^*_j+1$. 
By the definition of $i^*_j$, we have $\kappa_{i^*_j}(G) \geq j+1$ and 
$\kappa_{i^*+1}(G) \leq j$. Since $\kappa_{i^*_j} \geq j+1$,
$G$ is not $(j,i^*)$-colourable which means that $\lambda_j(G) \geq i^*+1$. 
Since $\kappa_{i^*_j+1}(G) \leq j$, $G$ is $(j,i^*+1)$-colourable which means
that $\lambda_j(G) \leq i^*+1$. Therefore we have $\lambda_j(G) = i^*+1$.
\qed

Theorem \ref{conjugate} allows us to convert $\hat{\kappa}(G)$ into 
$\hat{\lambda}(G)$ and vice versa. In particular, it implies that the two Ferrers 
diagrams of $\hat{\kappa}(G)$ and $\hat{\lambda}(G)$ have the same number of dots,
that is, 
\[\kappa_0(G)+\kappa_1(G)+\dots+\kappa_{\theta(G)-1}(G) = 
 \lambda_0(G)+\lambda_1(G)+\dots+\lambda_{\chi(G)-1}(G).\]
The number of vertices of the graph in Figure~\ref{g1} coincides with
the number of dots in either of the Ferrers diagrams in Figure~\ref{sqs_g1},
but this is not true in general. For instance, the 4-vertex graph $P_4$ has
$\hat{\kappa} = (2,1)$ whose Ferrers diagram has only three dots, whereas 
the 5-vertex graph $C_5$ has $\hat{\kappa} = (3,2,1)$ whose Ferrers diagram has 
six dots. However, we will show that for every cograph $G$, the number of 
vertice in $G$ is always equal to the number of dots in either of the Ferrers 
diagrams of $\hat{\kappa}(G)$ and $\hat{\lambda}(G)$ 
(see Theorem~\ref{T_CographsKappaSum_P}).

\begin{prop}\label{T_KappaLambdaUnionSpecific_P}
 Let $G,H$ be graphs and $k$ be a natural number. Then
 $$\lambda_k(G+H) = \lambda_k(G) + \lambda_k(H).$$
\end{prop}
\pf
  It suffices to show that $G+H$ admits a $(k,\lambda_k(G)+\lambda_k(H))$-colouring but not a $(k,\lambda_k(G)+\lambda_k(H)-1)$-colouring. For the first condition, let
  $$\seq{S}{k},\seq{C}{\lambda_k(G)}$$
  be a $(k,\lambda_k(G))$-colouring of $G$ and
  $$\seq{S'}{k},\seq{C'}{\lambda_k(H)}$$
  be a $(k,\lambda_k(H))$-colouring of $H$. Since $G+H$ contains no edges between $G$ and $H$, the sets $S_i \cup S'_{i}$ are independent for all $i$. Therefore
  $$S_1 \cup S'_1, \dots, S_k \cup S'_k, \seq{C}{\lambda_k(G)},\seq{C'}{\lambda_k(H)}$$
  is a $(k,\lambda_k(G) + \lambda_k(H))$-colouring of $G+H$.
  
  Now suppose, $G+H$ admits a $(k,\lambda_k(G)+\lambda_k(H)-1)$-colouring. Again, as $G+H$ has no edges between $G$ and $H$, every clique is completely contained in either $G$ or $H$. By the definition of $\lambda_k$, at least $\lambda_k(G)$ of the cliques are contained in $G$, while at least $\lambda_k(H)$ of the cliques are contained in $H$, implying that there are at least $\lambda_k(G) + \lambda_k(H)$ cliques in total, a contradiction. Therefore $G+H$ does not admit a $(k,\lambda_k(G)+\lambda_k(H)-1)$-colouring.
\qed

\begin{cor}\label{T_KappaLambdaUnionLambda_C}
 For any graphs $G,H$,
 $$\hat{\lambda}(G+H) = \hat{\lambda}(G) + \hat{\lambda}(H),$$
 where the addition is performed entrywise.
\qed
\end{cor}

The two sequqences $\hat{\lambda}(G)$ and $\hat{\lambda}(H)$ in 
Corollary~\ref{T_KappaLambdaUnionLambda_C} may have different length and if so
we append zeros to the shorter one to make them the same length.
By applying Proposition \ref{T_KappaLambdaUnionSpecific_P} and the fact that 
$\kappa_l(G) = \lambda_l(\overline{G})$, we can obtain the following equivalent statements for the join of two graphs.

\begin{prop}\label{T_KappaLambdaJoinSpecific_P}
 Let $G,H$ be graphs and $l$ a natural number. Then
 $$\kappa_l(G \vee H) = \kappa_l(G) + \kappa_l(H).$$
\end{prop}
\pf We have
  $$\kappa_l(G \vee H) = \lambda_l(\overline{G \vee H}) = \lambda_l(\overline{G} + \overline{H}) = \lambda_l(\overline{G}) + \lambda_l(\overline{H}) = \kappa_l(G) + \kappa_l(H).\vspace{-1.1cm}$$\vspace{0.4cm}
\qed

\begin{cor}\label{T_KappaLambdaJoinKappa_C}
 For any graphs $G,H$,
 $$\hat{\kappa}(G \vee H) = \hat{\kappa}(G) + \hat{\kappa}(H),$$
 where the addition is performed entrywise.
\qed
\end{cor}

Using the Ferrers diagrams, we can also compute $\hat{\lambda}(G \vee H)$ and 
$\hat{\kappa}(G+H)$. Consider $G \vee H$.  
By Corollary \ref{T_KappaLambdaJoinKappa_C}, $\hat{\kappa}(G \vee H)$ is obtained 
by adding $\hat{\kappa}(G)$ and $\hat{\kappa}(H)$ entrywise. 
In terms of the Ferrers diagram, we can picture this as putting 
the Ferrers diagrams of $\hat{\kappa}(G)$ and $\hat{\kappa}(H)$ beside each other 
with the rows lining up and moving all dots to the beginning of the row. 
This is equivalent to sorting the columns from largest to smallest. 
An example is given in Figure \ref{join}. We see that $\hat{\lambda}(G \vee H)$, 
being the conjugate of $\hat{\kappa}(G \vee H)$, is obtained by concatenating 
the sequence $\hat{\lambda}(G)$ with $\hat{\lambda}(H)$ and sorting the resulting 
sequence from largest to smallest. This way of concatenating two sequences
prompts the following definition of $*$.
\vspace{-4mm}

\begin{center}
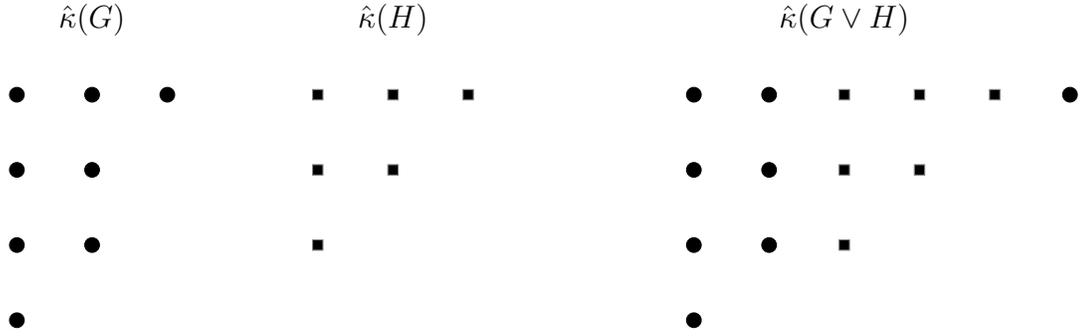
\begin{figure}[htb]
\center
\begin{tikzpicture}[>=latex]
\begin{pgfonlayer}{nodelayer}
                \node [style=blackvertex] (1) at (0,0){};
                \node [style=blackvertex] (2) at (1,0){};
                \node [style=blackvertex] (3) at (2,0){};
                \node [style=blackvertex] (4) at (0,-1){};
                \node [style=blackvertex] (5) at (1,-1){};
                \node [style=blackvertex] (6) at (0,-2){};
                \node [style=blackvertex] (7) at (1,-2){};
                \node [style=blackvertex] (8) at (0,-3){};

                \node [style=grayvertex] (9) at (4,0){};
                \node [style=grayvertex] (10) at (5,0){};
                \node [style=grayvertex] (11) at (6,0){};
                \node [style=grayvertex] (12) at (4,-1){};
                \node [style=grayvertex] (13) at (5,-1){};
                \node [style=grayvertex] (14) at (4,-2){};

                \node [style=blackvertex]  at (9,0){};
                \node [style=blackvertex] at (10,0){};
                \node [style=grayvertex] at (11,0){};
                \node [style=grayvertex] at (12,0){};
                \node [style=grayvertex] at (13,0){};
                \node [style=blackvertex] at (14,0){};

                \node [style=blackvertex] at (9,-1){};
                \node [style=blackvertex] at (10,-1){};
                \node [style=grayvertex] at (11,-1){};
                \node [style=grayvertex] at (12,-1){};

                \node [style=blackvertex] at (9,-2){};
                \node [style=blackvertex] at (10,-2){};
                \node [style=grayvertex] at (11,-2){};

                \node [style=blackvertex] at (9,-3){};

\node [style=textbox]  at (1, 1) {$\hat{\kappa}(G)$};
\node [style=textbox]  at (5, 1) {$\hat{\kappa}(H)$};
\node [style=textbox]  at (11, 1) {$\hat{\kappa}(G \vee H)$};

\end{pgfonlayer}
\end{tikzpicture}
\caption{\label{join}Ferrers diagrams of $\hat{\kappa}(G) = (3,2,2,1)$, $\hat{\kappa}(H) = (3,2,1)$ and $\hat{\kappa}(G \vee H) = (6,4,,3,1)$.}
\end{figure}
\end{center}

Let $\hat{a}$ and $\hat{b}$ be two non-increasing finite sequences of natural 
numbers. Then $\hat{a}*\hat{b}$ is the sequence obtained from concatenating 
$\hat{a}$ and $\hat{b}$ and sorting its entries from largest to smallest.

\begin{cor}\label{T_KappaLambdaJoinLambda_C}
 For any graphs $G,H$,
 $$\hat{\lambda}(G \vee H) = \hat{\lambda}(G) * \hat{\lambda}(H).$$
\qed
\end{cor}

\begin{cor}\label{T_KappaLambdaUnionKappa_C}
 For any graphs $G,H$,
 $$\hat{\kappa}(G+H) = \hat{\kappa}(G) * \hat{\kappa}(H).$$
\qed
\end{cor}

An interesting consequence of the above results is the following:

\begin{thm}\label{T_CographsKappaSum_P} 
 For any cograph $G$,
 $$\sum_{i \geq 0} \kappa_i(G) = \sum_{i \geq 0} \lambda_i(G) = \abs{V(G)}.$$
\end{thm}
\pf
  The first equality follows from Theorem \ref{conjugate}. We will show
  $$\sum_{i \geq 0} \kappa_i(G) = \abs{V(G)}$$
  by induction on the number of vertices of $G$. If $G = K_1$, then $\hat{\kappa}(G) = (1)$, thus the statement holds. Assume that $G \neq K_1$ and the statement holds for all cographs with fewer vertices than $G$. Suppose $G$ is disconnected. Then there exist cographs $G_1,G_2$ such that $G = G_1 + G_2$ and we have $\hat{\kappa}(G) = \hat{\kappa}(G_1) * \hat{\kappa}(G_2)$. By the definition of the operation $*$ and the induction hypothesis,
  \begin{align*}
   \sum_{i \geq 0} \kappa_i(G) &= \sum_{i \geq 0} \kappa_i(G_1) + \sum_{i \geq 0} \kappa_i(G_2)\\
                               &= \abs{V(G_1)} + \abs{V(G_2)}\\
                               &= \abs{V(G_1) \cup V(G_2)}\\
                               &= \abs{V(G)}.
  \end{align*}
  If on the other hand $G$ is connected, there exist cographs $G_1,G_2$ with $G = G_1 \vee G_2$ and we obtain as above
  \begin{align*}
   \sum_{i \geq 0} \kappa_i(G) &= \sum_{i \geq 0} (\kappa_i(G_1) + \kappa_i(G_2))\\
                              &= \sum_{i \geq 0} \kappa_i(G_1) + \sum_{i \geq 0} \kappa_i(G_2)\\
                              &= \abs{V(G)}.\qedhere
  \end{align*}
\qed

It is possible to draw a cograph $G$ on the Ferrers diagram of $\hat{\lambda}(G)$
by adding edges in such a way that the dots in each row form an independent set
and the dots in each column induce a clique in $G$. Such a drawing of $G$
is called a {\em Ferrers diagram representation} of $G$.

\vspace{-4mm}
\begin{center}
\begin{figure}[h!t]
\center
\begin{tikzpicture}[>=latex]
\begin{pgfonlayer}{nodelayer}
                \node [style=blackvertex] (1) at (0,0){};
                \node [style=blackvertex] (2) at (1,0){};
                \node [style=blackvertex] (3) at (2,0){};
                \node [style=blackvertex] (4) at (0,-1){};
                \node [style=blackvertex] (5) at (1,-1){};
                \node [style=blackvertex] (6) at (2,-1){};
                \node [style=blackvertex] (7) at (0,-2){};

                \node [style=grayvertex] (8) at (5,0){};
                \node [style=grayvertex] (9) at (6,0){};
                \node [style=grayvertex] (10) at (7,0){};
                \node [style=grayvertex] (11) at (5,-1){};
                \node [style=grayvertex] (12) at (6,-1){};
                \node [style=grayvertex] (13) at (5,-2){};
                \node [style=grayvertex] (14) at (5,-3){};

                \node [style=grayvertex]  at (10,0){};
                \node [style=blackvertex] at (11,0){};
                \node [style=blackvertex] at (12,0){};
                \node [style=blackvertex] at (13,0){};
                \node [style=grayvertex] at (14,0){};
                \node [style=grayvertex] at (15,0){};
                \node [style=grayvertex] at (10,-1){};
                \node [style=blackvertex] at (11,-1){};
                \node [style=blackvertex] at (12,-1){};
                \node [style=blackvertex] at (13,-1){};
                \node [style=grayvertex] at (14,-1){};
                \node [style=grayvertex] at (10,-2){};
                \node [style=blackvertex] at (11,-2){};
                \node [style=grayvertex] at (10,-3){};

\node [style=textbox]  at (1, 1) {$G_1$};
\node [style=textbox]  at (6, 1) {$G_2$};
\node [style=textbox]  at (12, 1) {$G_1 + G_2$};

\end{pgfonlayer}
\end{tikzpicture}
\caption{An illustration of the proof of Theorem \ref{T_CographsYoungDiagram_P}.}\label{F_CographsYoungUnion}
\end{figure}
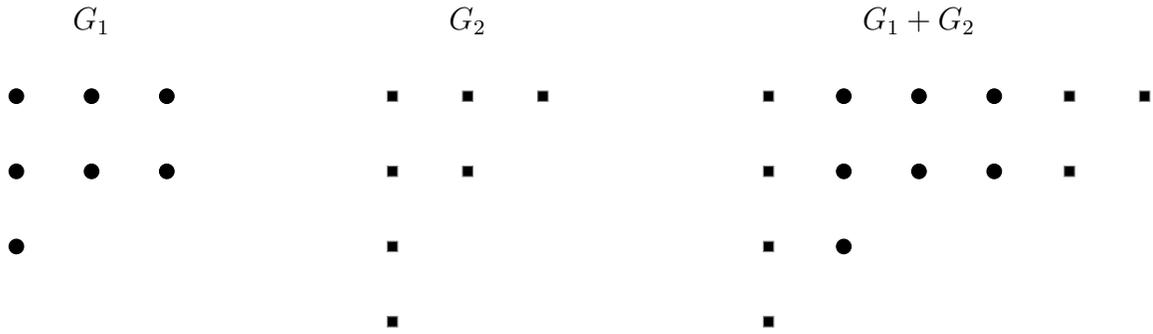
\end{center}

\begin{thm}\label{T_CographsYoungDiagram_P}
 Every cograph $G$ has a Ferrers diagram representation.
\end{thm}
\pf
  The proof is by induction on the number of vertices of $G$. If $G=K_1$, then the Ferrers diagram consists of a single point. Assume that $G$ has more than one vertex and every cograph on fewer vertices than $G$ has a Ferrers diagram representation. Suppose, $G$ is disconnected. Then there exist cographs $G_1,G_2$ with $G = G_1 + G_2$. By the induction hypothesis, $G_1$ and $G_2$ have a Ferrers diagram representation. Consider what happens if we write the two diagrams side by side (see the left side of Figure \ref{F_CographsYoungUnion} for an example). Each column completely belongs to either $G_1$ or $G_2$, therefore forms a clique. Each row may have vertices from both $G_1$ and $G_2$, but the vertices in each of the two graphs form an independent set, and therefore the whole row must form an independent set in $G_1+G_2$. The diagram we have might not be a Ferrers diagram, though. However, this can be remedied by permuting the columns. As this does not change the sets of rows and columns, we obtain a Ferrers diagram representation of $G$ (see the right side of Figure \ref{F_CographsYoungUnion} for the example).

  If $G$ is connected, there exist cographs $G_1,G_2$ with $G = G_1 \vee G_2$ and we obtain a Ferrers diagram representation of $G$ by writing the two Ferrers diagram representations of $G_1$ and $G_2$ on top of each other and permuting the rows instead.
\qed

\section{Cotrees}\label{S_CotreesBasicAlgorithms}

It follows from Theorem \ref{corneil} that every cograph $G$ on at least two 
vertices can be either written as $G = G_1 + G_2$ or $G = G_1 \vee G_2$ for some 
cographs $G_1, G_2$. These characteristic properties of cographs allow us to
represent a cograph $G$ as a tree, called the {\it cotree} of $G$,
cf. \cite{81CoLeSt}. The cotree of $G$, denoted by $T_G$, is a rooted tree where 
every internal node is labelled with either $0$ or $1$, which can be recursively 
constructed as follows.

\begin{itemize}
 \item[$\bullet$] If $G = K_1$, then we define $T_G$ to be the rooted tree on a single vertex.
 \item[$\bullet$] If $G$ is disconnected, let $\seq{G}{t}$ be the connected components of $G$. We take the cotrees of $G_1,G_2, \dots,G_t$ and add an edge from each of the roots to a new vertex, which we label with a $0$. The tree thus constructed, with the root at the new vertex, is the cotree $T_G$.
 \item[$\bullet$] If $G \neq K_1$ is connected, let $\overline{G_1},\overline{G_2},\dots,\overline{G_t}$ be the connected components of $\overline{G}$. We take the cotrees of $G_1,\dots,G_t$ and add an edge from each of the roots to a new vertex, which we label with a $1$. The tree thus constructed, with the root at the new vertex, is the cotree $T_G$.
\end{itemize}

We remark that the construction implies that the cotree $T_G$ for a cograph $G$
is unique. Every leaf of the cotree $T_G$ represents a vertex of $G$ and every internal node represents the subgraph of $G$ induced by the vertices that are descendents of that node. Every $0$-node represents a disconnected subgraph, every $1$-node a connected subgraph. By the construction, all children of a $0$-node represent connected cographs, thus are either $1$-nodes or leaves. Similarly the children of a $1$-node are either $0$-nodes or leaves. Also, we note that two vertices of the cograph are adjacent if and only if the lowest common ancestor of the corresponding leaves is a $1$-node. An example of a cotree is given in Figure \ref{F_CotreesCotree}. The corresponding cograph is shown in Figure \ref{F_CotreesCograph}, where the thick edges stand for complete adjacency.

\begin{conditional}{figures}{
 \begin{figure}[htb]
  \begin{center}
   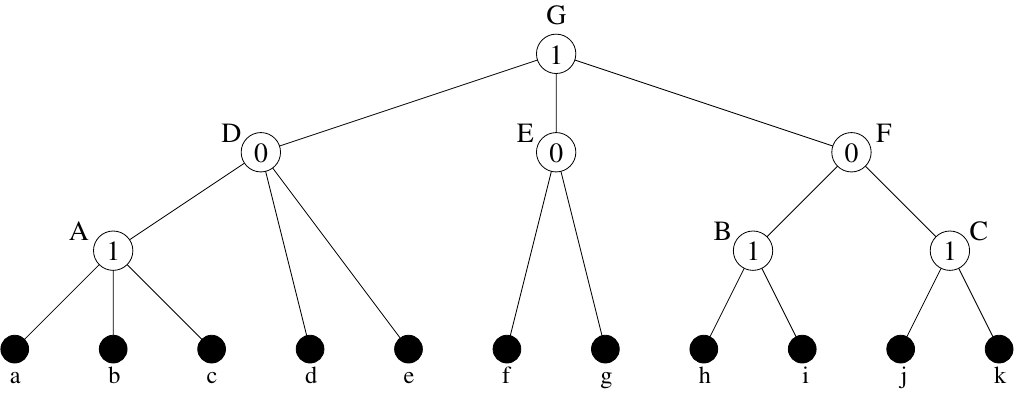
  \end{center}
  \caption{A cotree.}\label{F_CotreesCotree}
 \end{figure}
}\end{conditional}

\begin{conditional}{figures}{
 \begin{figure}[htb]
  \begin{center}
   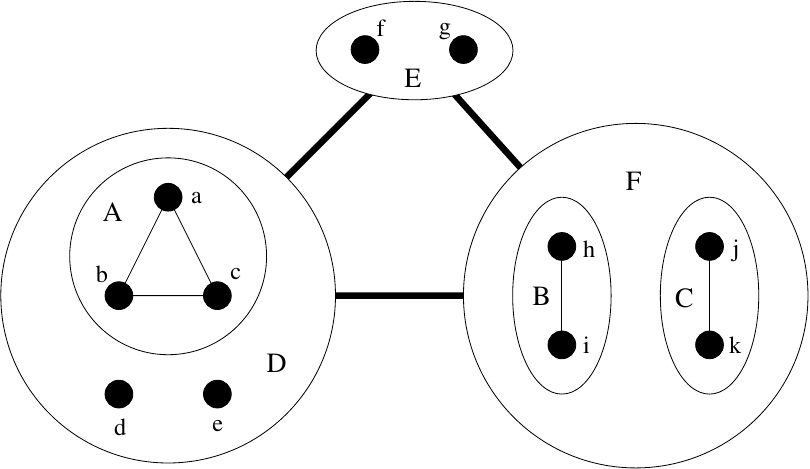
  \end{center}
  \caption{The cograph corresponding to the cotree from Figure \ref{F_CotreesCotree}.}\label{F_CotreesCograph}
 \end{figure}
}\end{conditional}

It is known that cotrees can be constructed in linear time (cf. \cite{85CoPeSt}).
Algorithms which are implemented on cotrees can be classified into two types,
depending on how they traverse on them. The bottom-up algorithm traverses 
the cotree from the leaves to the root, while the top-down algorithm traverses 
the cotree from the root to the leaves. 
Examples of bottom-up algorithms include calculating the chromatic number,
the cochromatic number, the number of cliques, and number of transitive 
orientations of a cograph (cf. \cite{81CoLeSt,94GiKrSt}).
We will give in Section \ref{S_CotreesColouring} two bottom-up algorithms, for 
calculating the sequence $\hat{\kappa}(G)$ and for constructing 
the Ferrers diagram representation of a cograph $G$ respectively.

Top-down algorithms are suited for example for finding an induced subgraph with 
certain properties, such as a maximum clique. Although not explicitly given
in \cite{81CoLeSt}, a top-down algorithm for finding a maximum clique in a cograph
can be derived from the bottom-up algorithm for the chromatic number of a cograph. We will also give in Section \ref{S_CotreesColouring} a top-down algorithm for 
finding an induced box cograph of a given dimension in a cograph.

\section{Calculating $\hat{\kappa}(G)$ and certificates}\label{S_CotreesColouring}

By Theorem~\ref{notejgt+}, a cograph is $(k,l)$-colourable if and only if it does
not contain a box cograph of dimension $k+1$ times $l+1$ as an induced subgraph.
In this section, we devise a bottom-up algorithm which calculates for any cograph 
$G$ the sequence $\hat{\kappa}(G)$ from which it can be determined whether $G$ is
$(k,l)$-colourable for any given $k, l$. In the case when $G$ is not
$(k,l)$-colourable a top-down algorithm will find an induced box cograph of
dimension $k+1$ times $l+1$ in $G$ which certifies that $G$ is not 
$(k,l)$-colourable.

Our bottom-up algorithm for the calculation of $\hat{\kappa}$ of a cograph relies
on the formulas for $\hat{\kappa}$ given 
in Corollaries \ref{T_KappaLambdaJoinKappa_C} and \ref{T_KappaLambdaUnionKappa_C}. 
The algorithm is similar in nature to the one presented in \cite{94GiKrSt} for 
the cochromatic number. However, the presentation is much simpler due to 
the formulas for $\hat{\kappa}$ established here.

\begin{alg}
 (KAPPA)\hfill
 \begin{itemize}
  \item[$\bullet$] INPUT: Cotree $T$ of a cograph $G$.
  \item[$\bullet$] INITIALIZATION: Assign $\hat{\kappa} = (1)$ to each leaf of $T$.
  \item[$\bullet$] $0$-NODE OPERATOR: $\hat{\kappa}(A) = \hat{\kappa}(A_1) * \hat{\kappa}(A_2) * \dots * \hat{\kappa}(A_t)$.
  \item[$\bullet$] $1$-NODE OPERATOR: $\hat{\kappa}(A) = \hat{\kappa}(A_1) + \hat{\kappa}(A_2) + \dots + \hat{\kappa}(A_t)$.
  \item[$\bullet$] OUTPUT: $\hat{\kappa}(G)$.
 \end{itemize}
\end{alg}
\pf [Proof of Correctness]
  The correctness follows directly from the fact that $\hat{\kappa}(K_1) = K_1$ and from the formulas for $\hat{\kappa}$ given in 
Corollaries \ref{T_KappaLambdaJoinKappa_C} and 
\ref{T_KappaLambdaUnionKappa_C} for the disjoint union and join of graphs.
\qed

An example of the output of KAPPA is shown in Figure \ref{F_CotreesColouringKL}, where the arguments of the leaves (all being $(1)$) have been omitted.

\begin{conditional}{figures}{
 \begin{figure}
  \begin{center}
   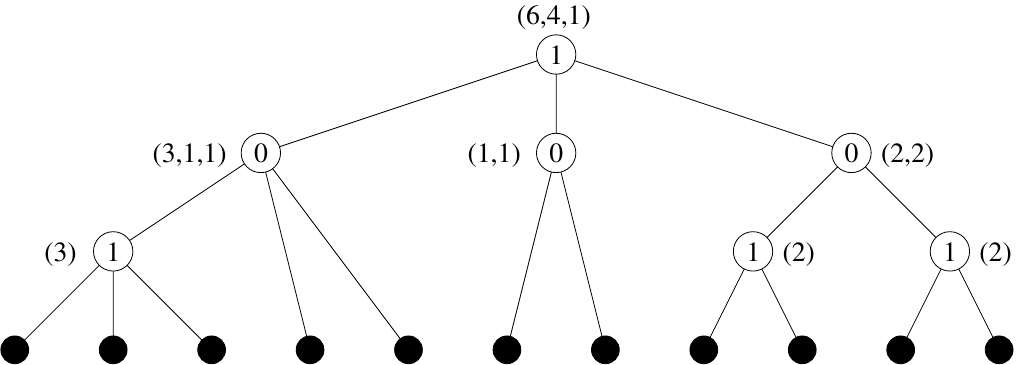
  \end{center}
  \caption{KAPPA for the cotree from Figure \ref{F_CotreesCotree}.}\label{F_CotreesColouringKL}
 \end{figure}
}\end{conditional}

As we can easily calculate the cochromatic number and bichromatic number from 
$\hat{\kappa}$, KAPPA can be seen as an algorithm for the $(k,l)$-colourability, 
the cochromatic number and the bichromatic number of a cograph. 

We will briefly discuss the complexity of KAPPA. Let $n$ be the number of vertices
of the cograph $G$. Then the number of operations ($*$ or $+$) performed by 
the algorithm is $O(n)$. We will show that each operation only needs time $O(n)$. 
Note that each sequence $\hat{\kappa}(A)$ has length $\abs{A}$. 
To calculate 
$\hat{\kappa}(A_1) * \hat{\kappa}(A_2)$, say, we need to sort the concatenated 
sequence of $\hat{\kappa}(A_1)$ and $\hat{\kappa}(A_2)$. Since both sequences are 
already sorted, we only need to scan each sequence once. As both sequences have 
length at most $n$, this can be done in $O(n)$. To calculate 
$\hat{\kappa}(A_1) + \hat{\kappa}(A_2)$, we need to perform at most $n$ additions,
which is also $O(n)$. Thus the algorithm can be implemented in time $O(n^2)$, 
which matches the complexity of the algorithm from \cite{94GiKrSt}. 
However, it is possible to slightly modify KAPPA to an algorithm that runs in time
$O(n \log n)$. For the sake of explanation, we consider a more general invariant
of cotrees, in which each internal node (again labeled either 0 or 1) has exactly 
two children. We call such a tree a {\em pseudocotree}. In general, a cograph 
can be represented by different pseudocotrees. But the number of nodes in
any pseudocotree for a cograph $G$ on $n$ vertices is $O(n)$. This means that 
the number of operations performed on any pseudocotree for $G$ is $O(n)$.
Suppose that $A$ is a node in the input pseudocotree for $G$ and $A_1$ and $A_2$ 
are the two children of $A$. We claim that $\hat{\kappa}(A)$ can be calculated in 
time $O(\min\set{\abs{A_1},\abs{A_2}})$. Indeed, assume without loss of generality
that $a = \abs{A_1} \leq \abs{A_2}$. If $A$ is a 1-node, then $\hat{\kappa}(A)$
can be obtained by adding the at most $a$ entries of $\hat{\kappa}(A_1)$ to the first
entries of $\hat{\kappa}(A_2)$, which can be done in time $O(a)$.  
Suppose that $A$ is a 0-node. The operation on $A$ requires to merge 
$\hat{\kappa}(A_1)$ into $\hat{\kappa}(A_2)$. To do it efficiently, we can
store $\hat{\kappa}(A_2)$ in the form of 
$(n_1^{\alpha_1}, n_2^{\alpha_2}, \dots, n_q^{\alpha_q})$ where 
$n_1 > n_2 > \cdots > n_q$. Thus it takes at most $O(a)$ scans of the entries
$n_i^{\alpha_i}$ for the merge of $\hat{\kappa}(A_1)$ into $\hat{\kappa}(A_2)$. 
Hence $\hat{\kappa}(A)$ can be calculated in time 
$O(\min\set{\abs{A_1},\abs{A_2}})$ for each node $A$ of the pseudocotree. 
This implies that the algorithm can be implemented to run in time 
$O(n \log n)$.

A similar algorithm as KAPPA can be devised to calculate $\hat{\lambda}(G)$
for a cograph $G$. This can be done simply by performing 
$\hat{\lambda}(A) = 
\hat{\lambda}(A_1) + \hat{\lambda}(A_2) + \dots + \hat{\lambda}(A_t)$ on each 
0-node $A$ and $\hat{\lambda}(A) = 
\hat{\lambda}(A_1) * \hat{\lambda}(A_2) * \dots * \hat{\lambda}(A_t)$ on each 
1-node $A$ of the cotree of $G$. The correctness of this algorithm is justified 
by Corollaries \ref{T_KappaLambdaUnionLambda_C} and 
\ref{T_KappaLambdaJoinLambda_C}.  

We can use KAPPA to establish a top-down algorithm that finds a certain 
box cograph. As a reminder, the obstructions for 
$(k-1,l-1)$-colourability of cographs are precisely the box cographs of dimension
$k$ times $l$, similar to the $k$-clique being the obstructions for 
$(k-1)$-colourability. We use $[r]^s$ to denote the sequence consisting of $s$ 
entries of $r$.

\begin{alg}
 (BOX COGRAPH)\hfill
 \begin{itemize}
  \item[$\bullet$] INPUT: Cotree $T$ of a cograph $G$ with $\kappa_{l-1}(G) \geq k$ and $\hat{\kappa}$ for each node of $T$.
  \item[$\bullet$] INITIALIZATION: Assign $c(G) = [k]^l$ to the root of $T$.
  \item[$\bullet$] $0$-NODE OPERATOR: For $c(A) = [r]^s$, set $c(A_i)$ such that
                   \begin{align*}
                    c(A_i) &= [r]^{s_i},\\
                    \kappa_{s_i-1}(G_i) &\geq r,\\
                    s_1 + s_2 + \dots + s_t = s.
                   \end{align*}
  \item[$\bullet$] $1$-NODE OPERATOR: For $c(A) = [r]^s$, set $c(A_i)$ such that
                   \begin{align*}
                    c(A_i) &= [r_i]^{s},\\
                    \kappa_{s-1}(G_i) &\geq r_i,\\
                    r_1 + r_2 + \dots + r_t = r.
                   \end{align*}
  \item[$\bullet$] OUTPUT: Leaves with $c=(1)$ inducing a box cograph 
       $H$ of dimension $k$ times $l$.
 \end{itemize}
\end{alg}
\pf [Proof of Correctness]
  We start by showing that the two operators are well-defined. Let $A$ be a node that got $c(A) = [r]^s$ assigned. By the initialization and the definition of the operators, we know that
  $$\kappa_{s-1}(A) \geq r.$$
  Suppose $A$ is a $0$-node. Then
  $$A = A_1 + A_2 + \dots + A_t$$
  and therefore
  $$\hat{\kappa}(A) = \hat{\kappa}(A_1) * \hat{\kappa}(A_2) * \dots * \hat{\kappa}(A_t).$$
  As $\kappa_{s-1}(A) \geq r$ implies that there are at least $s$ entries greater than or equal to $r$ in $\hat{\kappa}(A)$ and since $\hat{\kappa}(A)$ arises from the concatenation of $\hat{\kappa}(A_1), \hat{\kappa}(A_2),\dots, \hat{\kappa}(A_t)$, we know that we can find values $s_1,s_2,\dots,s_t$ such that $\hat{\kappa}(A_i)$ contains at least $s_i$ entries greater than or equal to $r$ and $s_1 + s_2 + \dots + s_t = s$. Therefore the $0$-operator is well-defined.

  If $A$ is a $1$-node, then
  $$A = A_1 \vee A_2 \vee \dots \vee A_t$$
  and we have
  $$\hat{\kappa}(A) = \hat{\kappa}(A_1) + \hat{\kappa}(A_2) + \dots + \hat{\kappa}(A_t),$$
  thus
  $$\kappa_{s-1}(A) = \kappa_{s-1}(A_1) + \kappa_{s-1}(A_2) + \dots + \kappa_{s-1}(A_t).$$
  Hence we can find values $r_1, r_2, \dots, r_t$ such that $\kappa_{s-1}(A_i) \geq r_i$ and $r_1 + r_2 + \dots + r_t = r$. Therefore the $1$-operator is well-defined.

  To show that the vertices with $c = (1)$ induce a box cograph of dimension
$k$ times $l$, we first note that for any $0$-node, the sum over the entries of $c(A)$ is
  $$rs = (r_1 + r_2 + \dots + r_t)s = r_1s + r_2s + \dots + r_ts,$$
  which is the sum over all entries of all $c(A_i)$. The same holds for $1$-nodes. Hence the sum over all entries of all sequences assigned to the leaves equals the sum over the entries of the sequence assigned to the root, which is $kl$. As the only possible assignments to the leaves are $(1)$ and the empty sequence $()$, we must have $kl$ leaves with $(1)$ assigned to them. It suffices to show that the graph $H$ induced by these leaves satisfies $\hat{\kappa}(H) = [k]^l$. To do so, we apply KAPPA to $T$, where we initialize the leaves by the arguments assigned to them by this algorithm. The operators of KAPPA are inverses to the ones of BOX COGRAPH. Therefore the output of KAPPA will be $c(G) = [k]^l$. It follows that
$\hat{\kappa}(H) = [k]^l$ and $H$ is a box cograph of dimension
$k$ times $l$.
\qed

An example of the output of BOX COGRAPH is shown in Figure \ref{F_CotreesColouringBoxCograph}, where $\hat{\kappa}$ is shown in round brackets to the left of each node and $c$ in square brackets to the right of each node. For the leaves, $\hat{\kappa}$ and $c$ have been omitted, except when $c = [1]$.

\begin{conditional}{figures}{
 \begin{figure}[ht]
  \begin{center}
   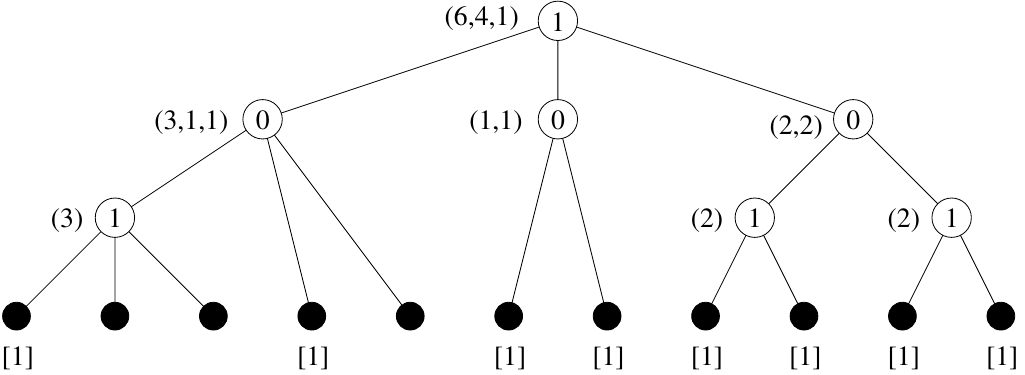
  \end{center}
  \caption{BOX COGRAPH with $k=4$ and $l=2$ for the cotree from Figure \ref{F_CotreesCotree}.}\label{F_CotreesColouringBoxCograph}
 \end{figure}
}\end{conditional}

As a final algorithm, we present a bottom-up algorithm, calculating 
the Ferrers diagram representation of a cograph.

\begin{alg}
 (FERRERS DIAGRAM)\hfill
 \begin{itemize}
  \item[$\bullet$] INPUT: Cotree $T$ of a cograph $G$.
  \item[$\bullet$] INITIALIZATION: Assign ${\cal Y} = \bullet$ to each leaf of $T$, labelled with the name of the leaf.
  \item[$\bullet$] $0$-NODE OPERATOR: ${\cal Y}(A)$ is the Ferrers diagram representation consisting of the columns of the Ferrers diagram representations ${\cal Y}(A_1),{\cal Y}(A_2),\dots,{\cal Y}(A_t)$, sorted by size.
  \item[$\bullet$] $1$-NODE OPERATOR: ${\cal Y}(A)$ is the Ferrers diagram representation consisting of the rows of the Ferrers diagram representations ${\cal Y}(A_1),{\cal Y}(A_2),\dots,{\cal Y}(A_t)$, sorted by size.
  \item[$\bullet$] OUTPUT: Ferrers diagram representation of $G$.
 \end{itemize}
\end{alg}
\pf [Proof of Correctness]
  The correctness follows directly from the proof of Theorem \ref{T_CographsYoungDiagram_P}.
\qed

An example of an output of FERRERS DIAGRAM is given in Figure \ref{F_CotreesColouringYoungDiagram}. The cotree is the one from Figure \ref{F_CotreesCotree} except that the labels of the internal nodes have been swapped for visual reasons (the graph represented by this cotree is the complement of the one from Figure \ref{F_CotreesCograph}).

\begin{conditional}{figures}{
 \begin{figure}
  \begin{center}
   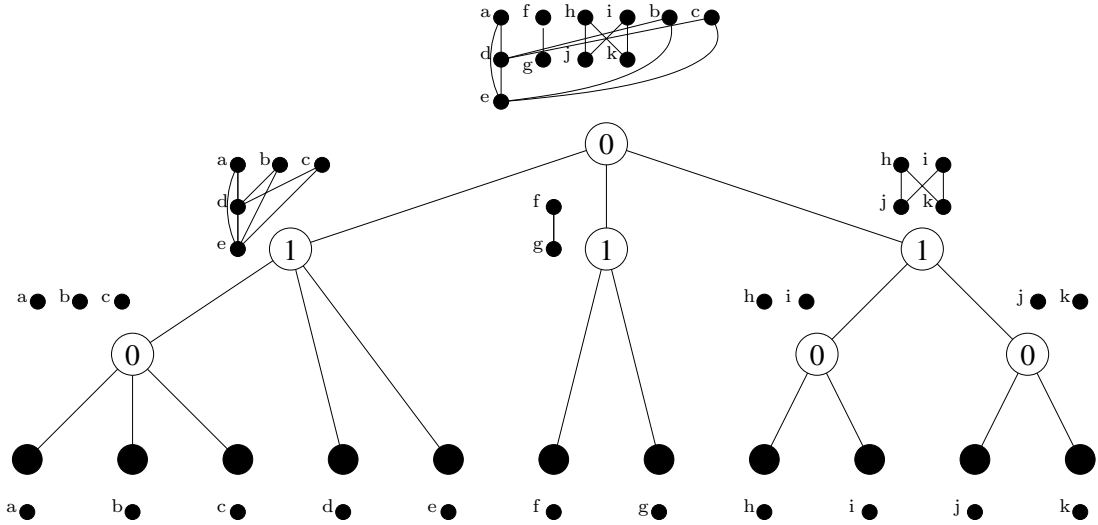
  \end{center}
  \caption{FERRERS DIAGRAM for the complement of the cotree from Figure \ref{F_CotreesCotree}.}\label{F_CotreesColouringYoungDiagram}
 \end{figure}
}\end{conditional}

Algorithm FERRERS DIAGRAM can be implemented in time $O(n\log n)$ similarly to
Algorithm KAPPA. We again describe the implementation on a pseudocotree. To store the array of vertices in the Ferrers diagram representation, we keep track of the right neighbour and bottom neighbour (in the representation not the graph) for each point. We also calculate $\hat{\kappa}$ and $\hat{\lambda}$ and the order of the induced subgraph at each node of the cotree. Suppose that $A$ is a node in the input pseudocotree and $A_1$ and $A_2$ are the two children of $A$. We claim that we can calculate the Ferrers diagram representation of $A$ from the Ferrers diagrams representation of $A_1$ and $A_2$ in $O(\min \{|A_1|,|A_2|\}$, implying that the algorithm can be implemented to run in time $O(n \log n)$. Assume without loss of generality $a = |A_1| \leq |A_2|$. If $A$ is a $0$-node we insert the columns of the Ferrers diagram representation of $A_1$ into the Ferrers diagram representation of $A_2$. To find the locations for inserting the columns requires $O(a)$ steps (using $\hat{\kappa}$). For each column we insert we need to update the right neighbours of its own vertices and of the vertices of the new column to its left, again $O(a)$ steps. Similarly, if $A$ is a $1$-node we can insert the rows of the Ferrers diagram representation of $A_1$ into that of $A_2$ in $O(a)$ steps, using $\hat{\lambda}$ to find the locations of insertion.

From the Ferrers diagram representation we can easily determine, whether a cograph is $(k,l)$-colourable and find an box cograph obstruction if not. For example, the graph in Figure \ref{F_CotreesColouringYoungDiagram} is not $(1,3)$-colourable. A box cograph of dimension
$2$ times $4$ is induced by the vertices $a,f,h,i,d,g,j,k$.

In this paper, we proved that the number of vertices in any cograph $G$ is equal
to the sum of entries in either of $\hat{\kappa}(G)$ or $\hat{\lambda}(G)$. 
We have seen examples (e.g., $P_4$ and $C_5$) for which this property does not
hold. It would be interesting to characterize all (perfect) graphs for which 
this property holds. 

{\bf Acknowledgement} We would to thank the anonymous referees for their 
helpful comments on the paper and suggestion of the above stated problem.

\bibliography{Ref}




\end{document}

%% file: Cotree.pdf_t
\begin{picture}(0,0)%
\includegraphics{Cotree.pdf}%
\end{picture}%
\setlength{\unitlength}{2072sp}%
\begingroup\makeatletter\ifx\SetFigFontNFSS\undefined%
\gdef\SetFigFontNFSS#1#2#3#4#5{%
  \reset@font\fontsize{#1}{#2pt}%
  \fontfamily{#3}\fontseries{#4}\fontshape{#5}%
  \selectfont}%
\fi\endgroup%
\begin{picture}(9270,3577)(766,-2258)
\end{picture}%

%% file: CotreeCograph.pdf_t
\begin{picture}(0,0)%
\includegraphics{CotreeCograph.pdf}%
\end{picture}%
\setlength{\unitlength}{2072sp}%
\begingroup\makeatletter\ifx\SetFigFontNFSS\undefined%
\gdef\SetFigFontNFSS#1#2#3#4#5{%
  \reset@font\fontsize{#1}{#2pt}%
  \fontfamily{#3}\fontseries{#4}\fontshape{#5}%
  \selectfont}%
\fi\endgroup%
\begin{picture}(7396,4281)(-187,-3444)
\end{picture}%

%% file: Cotree-k-l.pdf_t
\begin{picture}(0,0)%
\includegraphics{Cotree-k-l.pdf}%
\end{picture}%
\setlength{\unitlength}{2072sp}%
\begingroup\makeatletter\ifx\SetFigFontNFSS\undefined%
\gdef\SetFigFontNFSS#1#2#3#4#5{%
  \reset@font\fontsize{#1}{#2pt}%
  \fontfamily{#3}\fontseries{#4}\fontshape{#5}%
  \selectfont}%
\fi\endgroup%
\begin{picture}(9323,3324)(766,-1996)
\end{picture}%

%% file: Cotree-k-l-Obstruction.pdf_t
\begin{picture}(0,0)%
\includegraphics{Cotree-k-l-Obstruction.pdf}%
\end{picture}%
\setlength{\unitlength}{2072sp}%
\begingroup\makeatletter\ifx\SetFigFont\undefined%
\gdef\SetFigFont#1#2#3#4#5{%
  \reset@font\fontsize{#1}{#2pt}%
  \fontfamily{#3}\fontseries{#4}\fontshape{#5}%
  \selectfont}%
\fi\endgroup%
\begin{picture}(9353,3424)(751,-2398)
\put(6166,794){\makebox(0,0)[lb]{\smash{{\SetFigFont{8}{9.6}{\rmdefault}{\mddefault}{\updefault}{\color[rgb]{0,0,0}[4,4]}%
}}}}
\put(3421,-196){\makebox(0,0)[lb]{\smash{{\SetFigFont{8}{9.6}{\rmdefault}{\mddefault}{\updefault}{\color[rgb]{0,0,0}[1,1]}%
}}}}
\put(6121,-151){\makebox(0,0)[lb]{\smash{{\SetFigFont{8}{9.6}{\rmdefault}{\mddefault}{\updefault}{\color[rgb]{0,0,0}[1,1]}%
}}}}
\put(8821,-151){\makebox(0,0)[lb]{\smash{{\SetFigFont{8}{9.6}{\rmdefault}{\mddefault}{\updefault}{\color[rgb]{0,0,0}[2,2]}%
}}}}
\put(7921,-1051){\makebox(0,0)[lb]{\smash{{\SetFigFont{8}{9.6}{\rmdefault}{\mddefault}{\updefault}{\color[rgb]{0,0,0}[2]}%
}}}}
\put(9721,-1051){\makebox(0,0)[lb]{\smash{{\SetFigFont{8}{9.6}{\rmdefault}{\mddefault}{\updefault}{\color[rgb]{0,0,0}[2]}%
}}}}
\put(2071,-1051){\makebox(0,0)[lb]{\smash{{\SetFigFont{8}{9.6}{\rmdefault}{\mddefault}{\updefault}{\color[rgb]{0,0,0}[1]}%
}}}}
\end{picture}%

%% file: CotreeYoungDiagram.pdf_t
\begin{picture}(0,0)%
\includegraphics{CotreeYoungDiagram.pdf}%
\end{picture}%
\setlength{\unitlength}{2901sp}%
\begingroup\makeatletter\ifx\SetFigFont\undefined%
\gdef\SetFigFont#1#2#3#4#5{%
  \reset@font\fontsize{#1}{#2pt}%
  \fontfamily{#3}\fontseries{#4}\fontshape{#5}%
  \selectfont}%
\fi\endgroup%
\begin{picture}(9330,4421)(706,-2382)
\put(2521,299){\makebox(0,0)[lb]{\smash{{\SetFigFont{7}{8.4}{\rmdefault}{\mddefault}{\updefault}{\color[rgb]{0,0,0}d}%
}}}}
\put(2521,-61){\makebox(0,0)[lb]{\smash{{\SetFigFont{7}{8.4}{\rmdefault}{\mddefault}{\updefault}{\color[rgb]{0,0,0}e}%
}}}}
\put(2521,659){\makebox(0,0)[lb]{\smash{{\SetFigFont{7}{8.4}{\rmdefault}{\mddefault}{\updefault}{\color[rgb]{0,0,0}a}%
}}}}
\put(2881,659){\makebox(0,0)[lb]{\smash{{\SetFigFont{7}{8.4}{\rmdefault}{\mddefault}{\updefault}{\color[rgb]{0,0,0}b}%
}}}}
\put(3241,659){\makebox(0,0)[lb]{\smash{{\SetFigFont{7}{8.4}{\rmdefault}{\mddefault}{\updefault}{\color[rgb]{0,0,0}c}%
}}}}
\put(5221,299){\makebox(0,0)[lb]{\smash{{\SetFigFont{7}{8.4}{\rmdefault}{\mddefault}{\updefault}{\color[rgb]{0,0,0}f}%
}}}}
\put(5221,-61){\makebox(0,0)[lb]{\smash{{\SetFigFont{7}{8.4}{\rmdefault}{\mddefault}{\updefault}{\color[rgb]{0,0,0}g}%
}}}}
\put(7021,-511){\makebox(0,0)[lb]{\smash{{\SetFigFont{7}{8.4}{\rmdefault}{\mddefault}{\updefault}{\color[rgb]{0,0,0}h}%
}}}}
\put(7381,-511){\makebox(0,0)[lb]{\smash{{\SetFigFont{7}{8.4}{\rmdefault}{\mddefault}{\updefault}{\color[rgb]{0,0,0}i}%
}}}}
\put(9361,-511){\makebox(0,0)[lb]{\smash{{\SetFigFont{7}{8.4}{\rmdefault}{\mddefault}{\updefault}{\color[rgb]{0,0,0}j}%
}}}}
\put(9721,-511){\makebox(0,0)[lb]{\smash{{\SetFigFont{7}{8.4}{\rmdefault}{\mddefault}{\updefault}{\color[rgb]{0,0,0}k}%
}}}}
\put(8191,659){\makebox(0,0)[lb]{\smash{{\SetFigFont{7}{8.4}{\rmdefault}{\mddefault}{\updefault}{\color[rgb]{0,0,0}h}%
}}}}
\put(8551,659){\makebox(0,0)[lb]{\smash{{\SetFigFont{7}{8.4}{\rmdefault}{\mddefault}{\updefault}{\color[rgb]{0,0,0}i}%
}}}}
\put(8191,299){\makebox(0,0)[lb]{\smash{{\SetFigFont{7}{8.4}{\rmdefault}{\mddefault}{\updefault}{\color[rgb]{0,0,0}j}%
}}}}
\put(8551,299){\makebox(0,0)[lb]{\smash{{\SetFigFont{7}{8.4}{\rmdefault}{\mddefault}{\updefault}{\color[rgb]{0,0,0}k}%
}}}}
\put(4771,1919){\makebox(0,0)[lb]{\smash{{\SetFigFont{7}{8.4}{\rmdefault}{\mddefault}{\updefault}{\color[rgb]{0,0,0}a}%
}}}}
\put(4771,1559){\makebox(0,0)[lb]{\smash{{\SetFigFont{7}{8.4}{\rmdefault}{\mddefault}{\updefault}{\color[rgb]{0,0,0}d}%
}}}}
\put(4771,1199){\makebox(0,0)[lb]{\smash{{\SetFigFont{7}{8.4}{\rmdefault}{\mddefault}{\updefault}{\color[rgb]{0,0,0}e}%
}}}}
\put(6211,1919){\makebox(0,0)[lb]{\smash{{\SetFigFont{7}{8.4}{\rmdefault}{\mddefault}{\updefault}{\color[rgb]{0,0,0}b}%
}}}}
\put(6571,1919){\makebox(0,0)[lb]{\smash{{\SetFigFont{7}{8.4}{\rmdefault}{\mddefault}{\updefault}{\color[rgb]{0,0,0}c}%
}}}}
\put(5131,1919){\makebox(0,0)[lb]{\smash{{\SetFigFont{7}{8.4}{\rmdefault}{\mddefault}{\updefault}{\color[rgb]{0,0,0}f}%
}}}}
\put(5131,1469){\makebox(0,0)[lb]{\smash{{\SetFigFont{7}{8.4}{\rmdefault}{\mddefault}{\updefault}{\color[rgb]{0,0,0}g}%
}}}}
\put(5491,1919){\makebox(0,0)[lb]{\smash{{\SetFigFont{7}{8.4}{\rmdefault}{\mddefault}{\updefault}{\color[rgb]{0,0,0}h}%
}}}}
\put(5491,1559){\makebox(0,0)[lb]{\smash{{\SetFigFont{7}{8.4}{\rmdefault}{\mddefault}{\updefault}{\color[rgb]{0,0,0}j}%
}}}}
\put(5851,1919){\makebox(0,0)[lb]{\smash{{\SetFigFont{7}{8.4}{\rmdefault}{\mddefault}{\updefault}{\color[rgb]{0,0,0}i}%
}}}}
\put(5851,1559){\makebox(0,0)[lb]{\smash{{\SetFigFont{7}{8.4}{\rmdefault}{\mddefault}{\updefault}{\color[rgb]{0,0,0}k}%
}}}}
\put(811,-511){\makebox(0,0)[lb]{\smash{{\SetFigFont{7}{8.4}{\rmdefault}{\mddefault}{\updefault}{\color[rgb]{0,0,0}a}%
}}}}
\put(1171,-511){\makebox(0,0)[lb]{\smash{{\SetFigFont{7}{8.4}{\rmdefault}{\mddefault}{\updefault}{\color[rgb]{0,0,0}b}%
}}}}
\put(1531,-511){\makebox(0,0)[lb]{\smash{{\SetFigFont{7}{8.4}{\rmdefault}{\mddefault}{\updefault}{\color[rgb]{0,0,0}c}%
}}}}
\put(721,-2311){\makebox(0,0)[lb]{\smash{{\SetFigFont{7}{8.4}{\rmdefault}{\mddefault}{\updefault}{\color[rgb]{0,0,0}a}%
}}}}
\put(1621,-2311){\makebox(0,0)[lb]{\smash{{\SetFigFont{7}{8.4}{\rmdefault}{\mddefault}{\updefault}{\color[rgb]{0,0,0}b}%
}}}}
\put(2521,-2311){\makebox(0,0)[lb]{\smash{{\SetFigFont{7}{8.4}{\rmdefault}{\mddefault}{\updefault}{\color[rgb]{0,0,0}c}%
}}}}
\put(3421,-2311){\makebox(0,0)[lb]{\smash{{\SetFigFont{7}{8.4}{\rmdefault}{\mddefault}{\updefault}{\color[rgb]{0,0,0}d}%
}}}}
\put(4321,-2311){\makebox(0,0)[lb]{\smash{{\SetFigFont{7}{8.4}{\rmdefault}{\mddefault}{\updefault}{\color[rgb]{0,0,0}e}%
}}}}
\put(5221,-2311){\makebox(0,0)[lb]{\smash{{\SetFigFont{7}{8.4}{\rmdefault}{\mddefault}{\updefault}{\color[rgb]{0,0,0}f}%
}}}}
\put(6121,-2311){\makebox(0,0)[lb]{\smash{{\SetFigFont{7}{8.4}{\rmdefault}{\mddefault}{\updefault}{\color[rgb]{0,0,0}g}%
}}}}
\put(7021,-2311){\makebox(0,0)[lb]{\smash{{\SetFigFont{7}{8.4}{\rmdefault}{\mddefault}{\updefault}{\color[rgb]{0,0,0}h}%
}}}}
\put(7921,-2311){\makebox(0,0)[lb]{\smash{{\SetFigFont{7}{8.4}{\rmdefault}{\mddefault}{\updefault}{\color[rgb]{0,0,0}i}%
}}}}
\put(8821,-2311){\makebox(0,0)[lb]{\smash{{\SetFigFont{7}{8.4}{\rmdefault}{\mddefault}{\updefault}{\color[rgb]{0,0,0}j}%
}}}}
\put(9721,-2311){\makebox(0,0)[lb]{\smash{{\SetFigFont{7}{8.4}{\rmdefault}{\mddefault}{\updefault}{\color[rgb]{0,0,0}k}%
}}}}
\end{picture}%

%% file: dec2019.bbl
\providecommand{\bysame}{\leavevmode\hbox to3em{\hrulefill}\thinspace}
\providecommand{\MR}{\relax\ifhmode\unskip\space\fi MR }
\providecommand{\MRhref}[2]{%
  \href{http://www.ams.org/mathscinet-getitem?mr=#1}{#2}
}
\providecommand{\href}[2]{#2}
\begin{thebibliography}{10}

\bibitem{08AlSt}
N.~Alon and U.~Stav, \emph{The maximum edit distance from hereditary graph
  properties}, J. Combin. Theory Ser. B \textbf{98} (2008), no.~4, 672--697.

\bibitem{08AxKeMa}
M.~Axenovich, A.~{K\'ezdy}, and R.~Martin, \emph{On the editing distance of
  graphs}, J. Graph Theory \textbf{58} (2008), no.~2, 123--138.

\bibitem{97BoTh}
B.~{Bollob\'as} and A.~Thomason, \emph{Hereditary and monotone properties of
  graphs}, Algorithms Combin. \textbf{14} (1997), 70--78.

\bibitem{96Br}
A.~{Brandst\"adt}, \emph{Partitions of graphs into one or two independent sets
  and cliques}, Discrete Math. \textbf{152} (1996), no.~1-3, 47--54.

\bibitem{96Br_a}
A.~{Brandst{\"a}dt}, \emph{Corrigendum: {P}artitions of graphs into one or two
  independent sets and cliques}, Discrete Math. \textbf{186} (1998), no.~1-3,
  295.

\bibitem{81CoLeSt}
D.G. Corneil, H.~Lerchs, and L.K.~Stewart Burlingham, \emph{Complement
  reducible graphs}, Discrete Appl. Math \textbf{3} (1981), 163--174.

\bibitem{85CoPeSt}
D.G. Corneil, Y.~Perl, and L.K. Stewart, \emph{A linear recognition algorithm
  for cographs}, SIAM J. Comput. \textbf{14} (1985), 926--934.

\bibitem{05DeEkWe1}
M.~Demange, T.~Ekim, and D.~de~Werra, \emph{Partitioning cographs into cliques
  and stable sets}, Discrete Optim. \textbf{2} (2005), 145--153.

\bibitem{05DeEkWe2}
M.~{Demange}, T.~Ekim, and D.~de~Werra, \emph{{$(p,k)$-coloring} problems in
  line graphs}, Theoret. Comput. Sci. \textbf{349} (2005), no.~3, 462--474.

\bibitem{epple}
D.~Epple, \emph{Graph paritions and the bichromatic number}, PhD Thesis,
  University of Victoria, 2011.

\bibitem{10EpHu}
D.D.A. Epple and J.~Huang, \emph{A note on the bichromatic number of graphs},
  J. Graph Theory \textbf{65} (2010), no.~2, 263--269.

\bibitem{46ErSt}
P.~{Erd\H os} and A.H. Stone, \emph{On the structure of linear graphs}, Bull.
  Amer. Math. Soc. \textbf{52} (1946), 1087--1091.

\bibitem{07FeHeHo}
T.~Feder, P.~Hell, and W.~{Hochst\"attler}, \emph{Generalized colourings
  (matrix partitions) of cographs}, Graph Theory in Paris, Trends Math.,
  {Birkh\"auser}, Basel, 2007, pp.~149--167.

\bibitem{77FoHa}
S.~{F\"oldes} and P.L. Hammer, \emph{Split graphs having {D}ilworth number
  two}, Canad. J. Math. \textbf{29} (1977), no.~3, 666--672.

\bibitem{94GiKrSt}
J.~Gimbel, D.~Kratsch, and L.~Stewart, \emph{On cocolourings and cochromatic
  numbers of graphs}, Discrete Appl. Math. \textbf{48} (1994), 111--127.

\bibitem{04HeKlNoPr}
P.~Hell, S.~Klein, L.T. Nogueira, and F.~Protti, \emph{Partitioning chordal
  graphs into independent sets and cliques}, Discrete Appl. Math. \textbf{141}
  (2004), no.~1-3, 185--194.

\bibitem{77LeSt}
L.~Lesniak-Foster and H.J. Straight, \emph{The cochromatic number of a graph},
  Ars Combin. \textbf{3} (1977), 39--46.

\bibitem{93PrSt}
H.J. {Pr\"omel} and A.~Steger, \emph{Excluding induced subgraphs {II}: Extremal
  graphs}, Discrete Appl. Math. \textbf{44} (1993), 283--294.

\bibitem{turan}
P.~Tur\'an, \emph{On an extremal problem in graph theory}, Mat. Fiz. Lapok
  \textbf{48} (1941), 436--452 (in Hungarian).

\end{thebibliography}
